\documentclass[CEJM,PDF]{cej} 
\usepackage{layout}

\title{Weak Golbach's Conjecture from Isomorphic and Equivalent Odd Prime Number Functions}

\articletype{Research} 

\author{Luis A. Mateos\email{lamateos@amcomputersystems.com}
       }

\institute{AM Computer Systems Research, \\
           1070 Vienna, Austria
          }

\abstract{Mathematicians has been trying to prove the weak Goldbach's conjecture by adding prime numbers, as stated in the conjecture.
However, we believe that the solution does not need to be analytically solved. 
Instead of trying to add prime numbers to prove the conjecture, we 
developed a prime number function $P_{odd}(x)_{p>2}$, including $odd$ primes $p>2$, isomorphic and equivalent to a function $N_{odd}(x)_{n> 1}$, including $odd$ natural numbers greater than one, $n_{odd}> 1$, in which, the sum of three of its elements result in $odd$ numbers greater than 7, proving the conjecture.
}

\keywords{Prime Number \*\ Dynamical Prime Number function \*\ Goldbach's Conjecture \*\ Prime Sieves}
\msc{11P32, 37M10 }

\begin{document}
\maketitle
\section{Introduction }
The weak Goldbach's conjecture is a well know problem, "Every $odd$ number $n$ greater than 7 can be expressed as the sum of three $odd$ primes $p>2$" (An $odd$ prime $p>2$ may be used more than once in the same sum). 

\begin{equation}
\label{weakgold}
 \left( p_{p>2} + p_{p>2} + p_{p>2}\right)  \in n_{odd}> 7
\end{equation}

Mathematicians has been attacking the problem as stated, trying to add prime numbers.
Obtaining only partial results, proving the conjecture for sufficiently large $odd$ numbers \cite{hardy_cite} \cite{gove_cite}, in the range of $10^{32}$ \cite{bruno_cite}, if taking into account the Generalized Riemann Hypothesis (GRH).
Likewise, without taking into account the GRH, partial results has been obtained \cite{dimi_cite} \cite{boro_cite} \cite{chen_cite}, nevertheless, the sufficiently large $odd$ numbers increase drastically, in the range of $10^{7194}$ \cite{chen2_cite}.
Thus, still quantities not feasible to be checked by computers.

The Goldbach's conjecture is known to be true up to $10^{14}$. Deshouillers, te Riele and Saouter \cite{x1_cite} have checked it up to $10^{14}$ and Richstein \cite{x2_cite} up to $4 \times 10^{14}$.

In the same way, there has been parallel attempts to solve the conjecture, Tao \cite{tao_cite} in 2012 improved the result of Ramar\'{e} \cite{ramare_cite}, proving that every $odd$ number $N$ greater than 1 can be expressed as the sum of at most five primes, instead of six primes.

However, we believe the solution 
can be obtained in a different way.
By creating a mathematical function $P_{odd}(x)_{p>2}$ for all $odd$ primes, $p>2$, such that it is isomorphic (structurally identical) and equivalent to another mathematical function $N_{odd}(x)_{n> 1}$ for all $odd$ natural numbers greater than one, $n_{odd}> 1$, in which the sum of three of its elements will result in $odd$ numbers greater than 7, proving the conjecture.

\section{Prime Number Function $P(x)$}
Prime numbers are the building blocks of the positive integers, this was shown by Euler in his proof of the Euler product formula for the Riemann zeta function, Euler came up with a version of the sieve of Eratosthenes, better in the sense that each number was eliminated exactly once \cite{eulerproduct1_cite} \cite{eulerproduct2_cite}. 

\begin{equation}
\sum_{n=1}^{\infty} \frac{1}{n^s} = \prod_{p \, \mathrm{ prime}} \frac{1}{1-p^{-s}}
\end{equation}

The sieve of Eratosthenes is a simple algorithm for finding all prime numbers up to any given limit. It iteratively mark as composite the multiples of each prime, starting from that prime, while leaving without mark prime numbers \cite{sieves}.
In order to mimic this simple sieve algorithm as a mathematical function, the prime number function $p(x)$ includes the $sin(x)$ periodic function to indicate the periodicity of the prime multiples, $sin \left(\frac{1}{p}\right)$, starting from primes $p$ \cite{dynamicaleras}.

For a single prime number, the sieve equation is

\begin{equation}
\label{psfunctionsingle}
 p(x) = {p} + sin \left(\frac{1}{p}\right)
\end{equation}

While, for the entire set of primes, the equation can be represented as the set of all single prime functions $p(x)$

\begin{equation}
\label{psfunction}
 P(x) = \left\{p(x)\right\} =  \left\{{p} + sin \left(\frac{1}{p}\right)\right\}
\end{equation}

In order to visualize the prime number functions $P(x)$,  a one-dimensional number line is used, starting from zero and ending at infinite. Due to the nature of periodic functions, a zero-cross from the function over the number line, indicates the marking of a $composite$ number.
Consequently, by the inclusion of a periodic function in $p(x)$, the sieve works instantaneously from its starting prime $p$ to $\infty$. Likewise, by following a linear progression, starting at $p=2$ to $p\rightarrow\infty$, 
the multiples of primes are zero-cross, while prime number are decoded at each cycle of the function \cite{mateos1_cite}.

\begin{figure}[htp]
\caption{Zero-cross numbers by P(x).
In red $p(x)$ for $p=2$, in black $P(x)$ for $p>2$ ($odd$ primes);
white squares zero-cross represent $composite$ numbers; white squares zero-cross by only $p(x)$ for $p=2$ (red color) represent numbers power of two $2^m$, for $m>1$. \label{fig_progression}}
\includegraphics[width=1\textwidth]{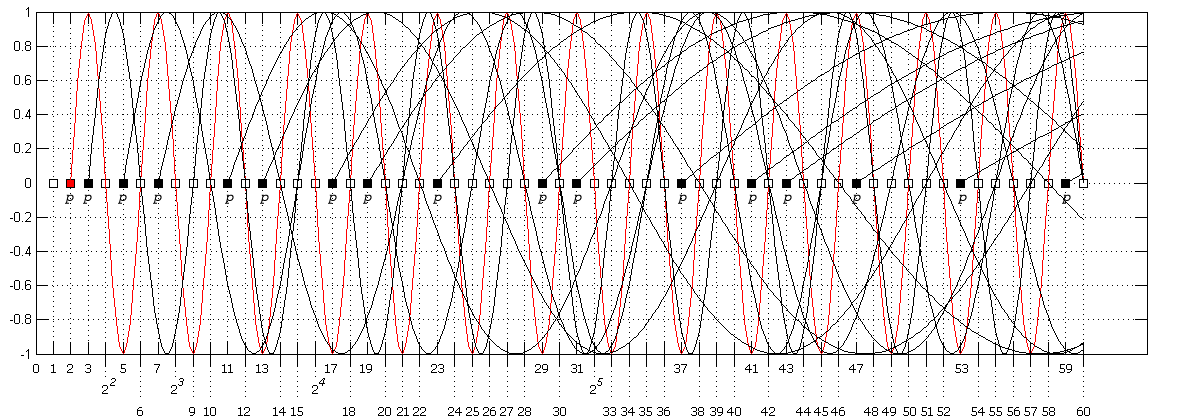}
\end{figure}

At each number $n>1$, the function $P(x)$ evaluates if the number $n$ in consideration has been zero-cross or not. If a number $n$ has not been zero-cross, then the number is a prime $p$, this process we called $prime$ $decoding$. 
While, if a number $n$ has been zero-cross, then the number is marked as $composite$, meaning that a previous prime or primes $p$ has been decoded and this zero-cross is a multiple of such prime or primes $p(x)$, this process we called $composite$ $encoding$.
Once a prime number $p$ is decoded, the function $p(x)$ starts the expansion of the prime into its multiples, mimicking the sieve of Eratosthenes, zero-crossing $composite$ numbers and decoding primes by leaving intact the number line at prime numbers $p$, as shown in figure \ref{fig_progression}. 

\begin{figure}[htp]
\caption{$P(x)=\left\{2+ sin \left(\frac{1}{2}\right) \right\}$ starts at $p=2$ and zero-cross multiples of $2$. In red $p=2$. \label{fig_2}}
\includegraphics[width=1\textwidth]{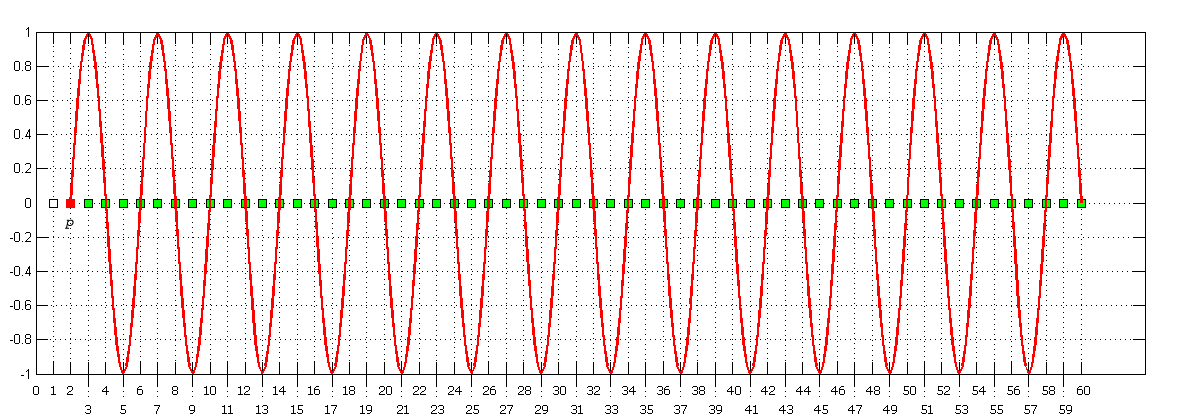}
\end{figure}

\begin{figure}[htp]
\caption{$P(x)=\left\{2+ sin \left(\frac{1}{2}\right), 3+ sin \left(\frac{1}{3}\right) \right\}$ zero-cross multiples of $2$ and $3$. In red $p=3$. \label{fig_3}}
\includegraphics[width=1\textwidth]{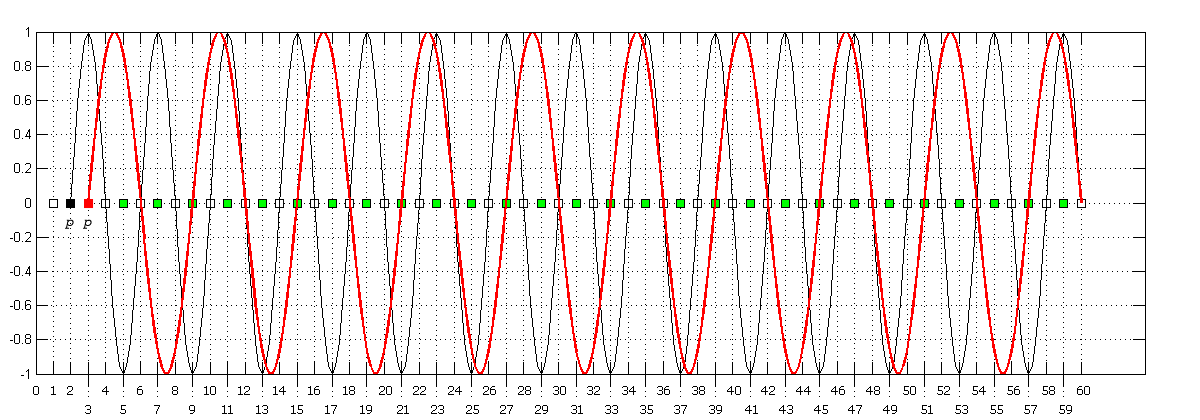}
\end{figure}

\begin{figure}[htp]
\caption{$P(x)= \left\{ 2+ sin \left(\frac{1}{2}\right), 3+ sin \left(\frac{1}{3}\right),  5+ sin \left(\frac{1}{5}\right) \right\}$ zero-cross multiples of $2$, $3$ and $5$. In red $p=5$. \label{fig_5}}
\includegraphics[width=1\textwidth]{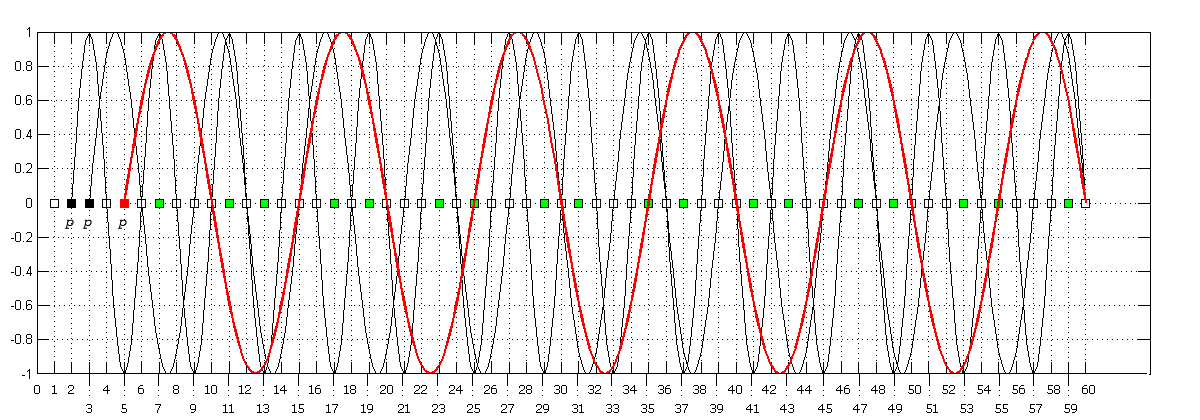}
\end{figure}

\begin{figure}[htp]
\caption{$P(x)= \left\{ 2+ sin \left(\frac{1}{2}\right), 3+ sin \left(\frac{1}{3}\right), 5+ sin \left(\frac{1}{5}\right), 7+ sin \left(\frac{1}{7}\right) \right\}$ zero-cross multiples of $2$, $3$, $5$ and $7$. In red $p=7$. \label{fig_7}}
\includegraphics[width=1\textwidth]{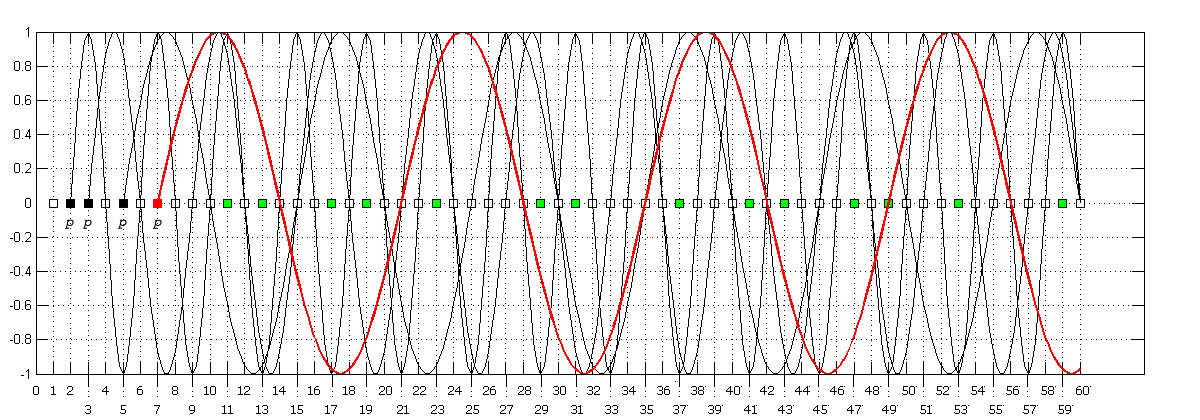}
\end{figure}

\begin{figure}[htp]
\caption{$P(x)= \left\{ 2+ sin \left(\frac{1}{2}\right), 3+ sin \left(\frac{1}{3}\right), 5+ sin \left(\frac{1}{5}\right), 7+ sin \left(\frac{1}{7}\right), 11+ sin \left(\frac{1}{11}\right) \right\}$ zero-cross multiples of $2$, $3$, $5$, $7$ and $11$. In red $p=11$. \label{fig_11}}
\includegraphics[width=1\textwidth]{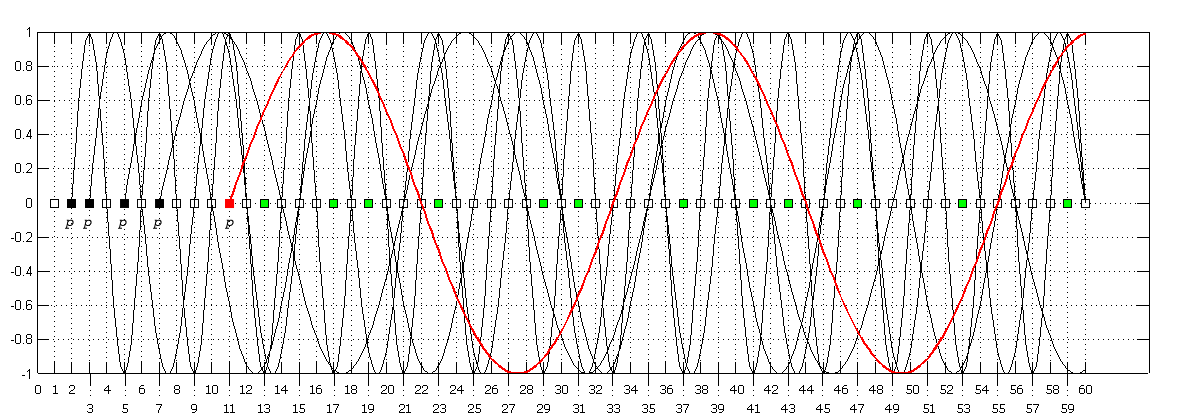}
\end{figure}

The process is as follows:
Starting at number $p=2$ the function $p(x)$$=2+ sin \left(\frac{1}{2}\right)$ zero-crosses all multiples of $2$, while the maximum value of the function exist at all $odd$ numbers, $n_{odd}>1$, as shown in figure \ref{fig_2}.
At number $p=3$ no zero-cross occur. Thus, the function $p(x)$$=3+ sin \left(\frac{1}{3}\right)$ starts zero-crossing $odd$ and $even$ multiples of $3$, as shown in figure \ref{fig_3}.
At numbers $n=4$ $n=6$,  $n=8$, $n=10$ and all $even$ numbers $n_{even}>2$ a zero-cross occurs from the function $p(x)$$=2 + sin \left(\frac{1}{2}\right)$ meaning $composite$ numbers, fig. \ref{fig_2}.
At number  $p=5$ no zero-cross occur. Therefore, the function $p(x)$$=5+ sin \left(\frac{1}{5}\right)$ starts to zero-cross $odd$ and $even$ multiples of $5$, as shown in figure \ref{fig_5}.
At number $p=7$ no zero-cross occur. Thus, the function $p(x)$$=7+ sin \left(\frac{1}{7}\right)$ starts to zero-cross $odd$ and $even$ multiples of $7$, as shown in figure \ref{fig_7}.
At number $n=9$ a zero-cross occurs from the function $p(x)$$=3+ sin \left(\frac{1}{3}\right)$ meaning a $composite$ number has been found, fig. \ref{fig_3}.
At number $p=11$ no zero-cross occur. Thus the function $p(x)$$=11+ sin \left(\frac{1}{11}\right)$ starts to zero-cross $odd$ and $even$ multiples of $11$, as shown in figure \ref{fig_11}. And so on.

Consequently, all $p(x)$ prime functions for $p>2$, zero-cross $odd$ and $even$ multiples of $p>2$.

\subsection{Natural Number Function $N(x)_{n>1}$}
If, instead of include only prime numbers $p$ in the prime number function $P(x)$, natural numbers greater than one are included, $n> 1$. 
The function changes to $N(x)_{n>1}$, starting from natural numbers greater than one and zero-crossing the same numbers as the prime number function $P(x)$, the set of $composite$ numbers. 

The sieve equation for a single natural number function is

\begin{equation}
\label{singlenatural}
 n(x)_{n>1} = {n_{n>1}} + sin \left(\frac{1}{n_{n>1}}\right)
\end{equation}

While, for the entire set of natural number functions $n(x)_{n>1}$, the equation can be represented as

\begin{equation}
\label{natural}
 N(x)_{n>1} = \left\{n(x)_{n>1}\right\} =  \left\{{n_{n>1}} + sin \left(\frac{1}{n_{n>1}}\right)\right\}
\end{equation}

There are two possible cases when applying this function:

{\bf Case I}.
If the function $N(x)_{n>1}$ follows the sieve rules, meaning that 
single functions $n(x)_{n>1}$ only start from non-zero-cross numbers greater than one (primes), then both functions are equal

\begin{equation}
\label{naturalsieve}
N(x)_{n>1}=P(x)
\end{equation}

{\bf Case II}.
On the other hand, if we change the rule, meaning that there exist a function $n(x)_{n>1}$ at every point greater than one over the number line, regardless if the point in consideration has been zero-cross or not, then the functions are equivalent and isomorphic, but not equal.



\begin{equation}
\label{isomorphic1}
 N(x)_{n>1} = \left\{n(x)_{n>1}\right\} =  \left\{{n_{n>1}} + sin \left(\frac{1}{n_{n>1}}\right)\right\} \iff  P(x) = \left\{p(x)\right\} =  \left\{{p} + sin \left(\frac{1}{p}\right)\right\}
\end{equation}

However, in both cases, the sieve functions, $N(x)_{n>1}$ and $P(x)$,  will lead to the same result, the decoding of prime numbers while encoding $composite$ numbers. 

\subsection{Odd Natural Number Function $N_{odd}(x)_{n>1}$}
The weak Goldbach's conjecture only takes into account $odd$ numbers, therefore, all the involved function must have as input and output $odd$ numbers greater than one.
For this purpose
the function $N(x)_{n>1}$ is modified to $N_{odd}(x)_{n>1}$, starting from $odd$ natural numbers greater than one and zero-crossing only its $odd$ multiples 
by doubling its period, from $odd$ periods $sin \left(\frac{1}{n_{odd_{n>1}}}\right)$ to $even$ periods $sin \left(\frac{1}{2n_{odd_{n>1}}}\right)$.

The sieve equation for a single $odd$ natural number function is

\begin{equation}
\label{oddnumberfunction}
 n_{odd}(x)_{n> 1} = {n_{odd_{n>1}}} + sin \left(\frac{1}{2n_{odd_{n>1}}}\right)
\end{equation}

While, for the entire set of $odd$ natural number functions, the equation can be represented as

\begin{equation}
\label{oddnumberfunctionall}
 N_{odd}(x)_{n>1} = \left\{n_{odd}(x)_{n> 1}\right\} =  \left\{{n_{odd_{n>1}}} + sin \left(\frac{1}{2n_{odd_{n>1}}}\right)\right\}
\end{equation}

Similar to the function $N(x)_{n>1}$, there are two ways to apply this developed function.

{\bf Case I}. If the function $N_{odd}(x)_{n>1}$ follows the sieve rules, in which, single $odd$ functions $n_{odd}(x)_{n> 1}$ only start from non-zero-cross $odd$ numbers ($odd$ primes).

{\bf Case II}. If the function $N_{odd}(x)_{n>1}$ does not follows the sieve rules, meaning that there exist a function $n_{odd}(x)_{n> 1}$ at every $odd$ point greater than one over the number line, regardless, if the point in consideration has been zero-cross or not.

However, if the function $N_{odd}(x)_{n>1}$ follows the sieve rules or not, the result will be the same, the decoding of $odd$ prime numbers, $p>2$, while zero-crossing $odd$ $composite$ numbers.

\section{Weak Golbach's Conjecture from Odd Prime Functions $P_{odd}(x)_{p>2}$}
The weak Goldbach's conjecture states: "Every $odd$ number $n$ greater than 7 can be expressed as the sum of three $odd$ primes $p>2$". (An $odd$ prime $p>2$ may be used more than once in the same sum).

If instead of $odd$ primes $p>2$, the conjecture states $odd$ numbers greater than one, $n_{odd}>1$.

"Every $odd$ number $n$ greater than 7 can be expressed as the sum of three $odd$ numbers greater than one, $n_{odd}>1$". (An $odd$ number greater than one, $n_{odd}>1$ may be used more than once in the same sum).
Then this modified conjecture is proven.

\begin{corollary}
\label{weakgoldn}

The sum of three $odd$ numbers bigger than one, $n_{odd}>1$, will result in $odd$ numbers greater than 7, $n_{odd}>7$

\begin{equation}
\label{weakgold}
 \left( n_{odd}> 1 + n_{odd}> 1 + n_{odd}> 1\right) \in n_{odd} > 7
\end{equation}

\end{corollary}

\begin{proposition}\label{Standard-stuff}

The mathematical function $P_{odd}(x)_{p>2}$ for all $odd$ primes, $p>2$, is isomorphic (structurally identical) and equivalent to the mathematical function $N_{odd}(x)_{n> 1}$ for all $odd$ natural numbers greater than one, $n_{odd}> 1$.

\begin{equation}
\label{isomorph}
 P_{odd}(x)_{p>2} \iff N_{odd}(x)_{n> 1}
\end{equation}

\end{proposition}

\begin{figure}[htp]
\caption{$P_{odd}(x)_{p > 2} = \left\{ {3} + sin \left(\frac{1}{6}\right)\right\}$ starts at $p=3$ and zero-cross $odd$ multiples of $3$. In blue $p=3$. \label{fig_podd3}}
\includegraphics[width=1\textwidth]{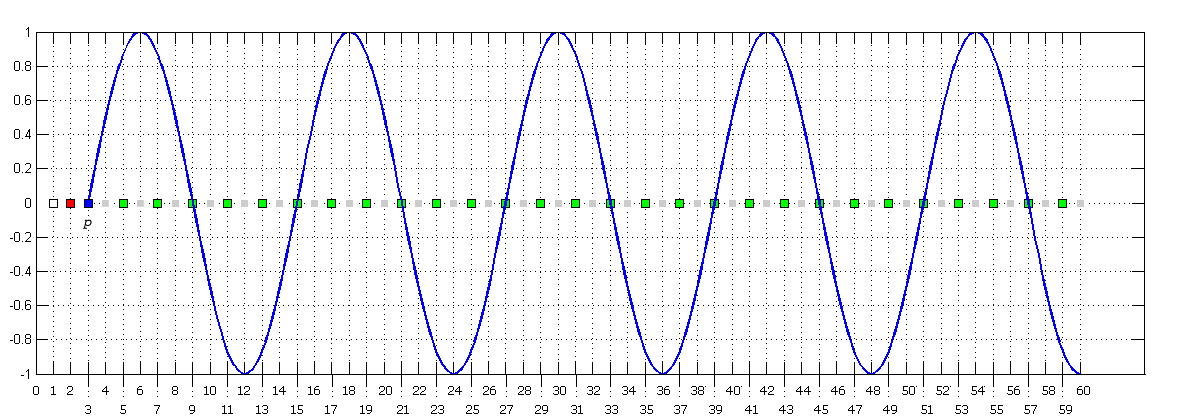}
\end{figure}

\begin{figure}[htp]
\caption{$P_{odd}(x)_{p > 2} = \left\{ {3} + sin \left(\frac{1}{6}\right) , {5} + sin \left(\frac{1}{10}\right) \right\}$ zero-cross $odd$ multiples of $p=3$ and $p=5$. In black $p=3$, in blue $p=5$. \label{fig_podd35}}
\includegraphics[width=1\textwidth]{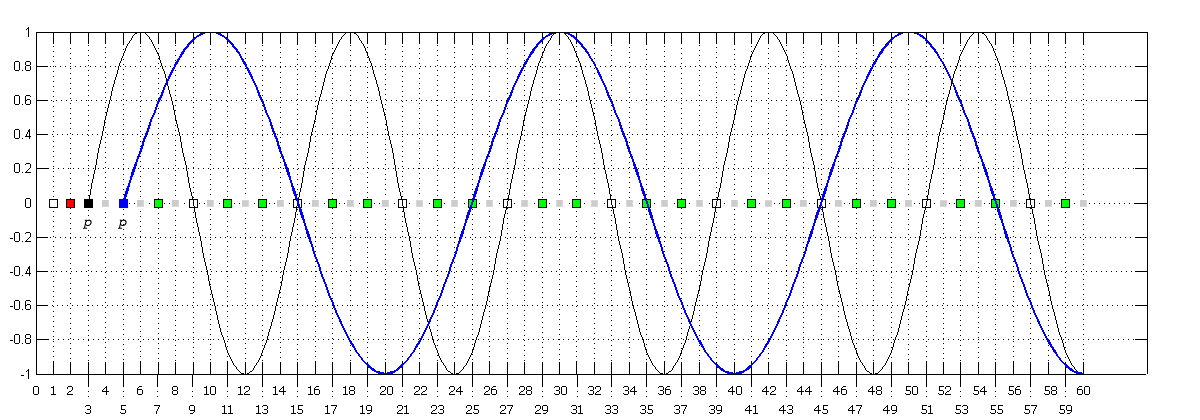}
\end{figure}

\begin{figure}[htp]
\caption{$P_{odd}(x)_{p > 2} = \left\{ {3} + sin \left(\frac{1}{6}\right) , {5} + sin \left(\frac{1}{10}\right) , {7} + sin \left(\frac{1}{14}\right) \right\}$ zero-cross $odd$ multiples of $p=3$, $p=5$ and $p=7$. In black $p=3$ and $p=5$, in blue $p=7$. \label{fig_podd357}}
\includegraphics[width=1\textwidth]{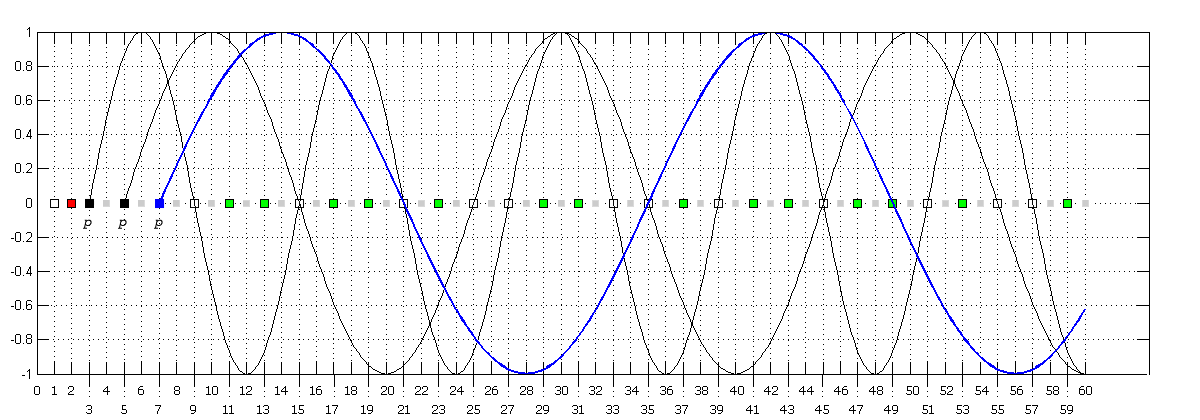}
\end{figure}

\begin{figure}[htp]
\caption{$P_{odd}(x)_{p > 2} = \left\{ {3} + sin \left(\frac{1}{6}\right) , {5} + sin \left(\frac{1}{10}\right) , {7} + sin \left(\frac{1}{14}\right) , {11} + sin \left(\frac{1}{22}\right) \right\}$ zero-cross $odd$ multiples of $p=3$, $p=5$, $p=7$ and $p=11$. In black $p=3$,  $p=5$,  $p=7$, in blue $p=11$. \label{fig_podd35711}}
\includegraphics[width=1\textwidth]{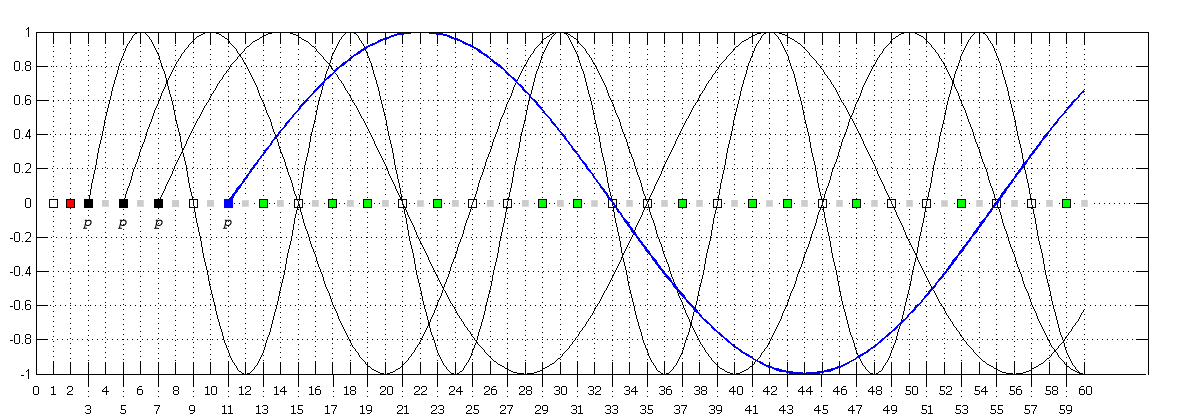}
\end{figure}

\begin{figure}[htp]
\caption{$P_{odd}(x)_{p > 2}$ zero-crosses $odd$ multiples of $p>2$, starting from $odd$ primes $p>2$. In red $p=2$, in black $p>2$. \label{fig_podds}}
\includegraphics[width=1\textwidth]{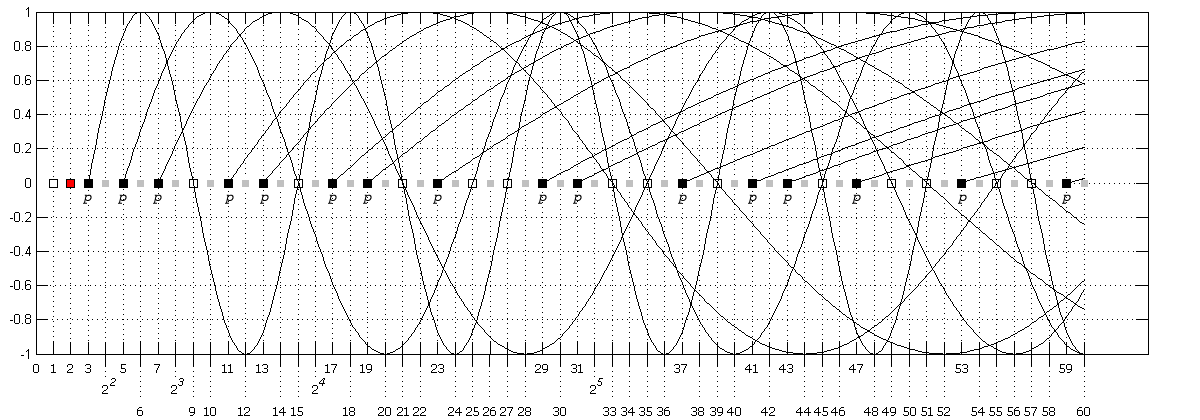}
\end{figure}

\begin{proof}
By definition, to obtain a prime number involves a process. "A prime number is a natural number greater than 1 that has no positive divisors other than 1 and itself".

Observe that the process to obtain primes $p$ can be written as the function $P(x)$, equation \ref{psfunction}, similar to the sieve of Eratosthenes, which iteratively mark as $composite$ the multiples of each prime, starting from that prime, while leaving without mark prime numbers \cite{sieves} \cite{dynamicaleras}.

The prime number function $P(x)$, takes into account all $even$ and $odd$ primes. If taking out of consideration $p(x)$, when $p=2$, then the function only includes $odd$ primes $P(x)_{p>2}$, and the zero-cross $composite$ numbers change, excluding only powers of two $\sum_{m>1}^{\infty} 2^{m}$, as shown in figure \ref{fig_progression}.

In order to create an $odd$ prime function $P_{odd}(x)_{p > 2}$ isomorphic and equivalent to $N_{odd}(x)_{n>1}$, equation \ref{oddnumberfunctionall}. 
It is necessary to modify the function $P(x)_{p>2}$, to zero-cross only $odd$ $composite$ numbers. This is done by doubling the period of the function $P(x)_{p>2}$, so the function instead of $odd$ periods $ \left(\frac{1}{p_{p>2}}\right)$, consist of $even$ periods $ \left(\frac{1}{2p_{p>2}}\right)$. Consequently, the function for a single $odd$ prime number is modified to

\begin{equation}
\label{psfunctionsingledd}
 p_{odd}(x)_{p > 2} = {p_{p>2}} + sin \left(\frac{1}{2p_{p>2}}\right)
\end{equation}

Starting at each $odd$ prime number, $p>2$, and zero-crossing only its $odd$ multiples. While the maximum value of the function exists at $even$ multiples of $odd$ primes, as shown in figure \ref{fig_podds}. 
In the same way, the function $P_{odd}(x)_{p > 2}$ represents the set of all single $odd$ prime functions $p_{odd}(x)_{p > 2}$

\begin{equation}
\label{eq-psfunction2}
 P_{odd}(x)_{p > 2} = \left\{ p_{odd}(x)_{p > 2} \right\} =   \left\{ {p_{p>2}} + sin \left(\frac{1}{2p_{p>2}}\right) \right\}
\end{equation}

Starting at $p=3$, the function $p_{odd}(x)= \left[ {3} + sin \left(\frac{1}{6}\right) \right]$ zero-crosses $odd$ multiples of $3$, while the maximum value of the function exist at $even$ multiples of $3$, as shown in figure \ref{fig_podd3}.
At number $p=5$ no zero-cross occur. Thus, the function $p_{odd}(x) = \left[ {5} + sin \left(\frac{1}{10}\right) \right]$ starts to zero-cross $odd$ multiples of $5$, as shown in figure \ref{fig_podd35}.
At number $p=7$ no zero-cross occur. Thus, the function $p_{odd}(x) = \left[ {7} + sin \left(\frac{1}{14}\right) \right]$ starts to zero-cross $odd$ multiples of $7$, as shown in figure \ref{fig_podd357}.
At number $n=9$ a zero-cross occurs from the function $p_{odd}(x)$ for $p= 3$, meaning the number is $composite$, fig. \ref{fig_podd3}. 
At number $p=11$ no zero-cross occur. Therefore, the function $p_{odd}(x) = \left[ {11} + sin \left(\frac{1}{22}\right) \right]$ starts to zero-cross $odd$ multiples of $11,$ as shown in figure \ref{fig_podd35711}. And so on.

Consequently, the function $P_{odd}(x)_{p > 2}$, zero-crosses the entire set of $odd$ composite numbers, starting from $odd$ primes, as shown in figure \ref{fig_podds}. 

Thus, 
the function $P_{odd}(x)_{p > 2}$ is equal to $N_{odd}(x)_{n>1}$, for {\bf case I}, when the sieve rules are applied, where single natural number functions $n_{odd}(x)_{n > 1}$ only start from non-zero-cross $odd$ numbers ($odd$ primes)

\begin{equation}
\label{oddequal}
 P_{odd}(x)_{p > 2} = N_{odd}(x)_{n>1}
\end{equation}

While, for {\bf case II}, there exist a function $n_{odd}(x)_{n>1}$ at every $odd$ point greater than one over the number line, regardless, if the point in consideration has been zero-cross or not.
Then the functions are isomorphic and equivalent, equation \ref{isomorph}.


\end{proof}






\begin{proposition}\label{Standard-stuff2}
The mathematical function $P_{odd}(x)_{p>2}$ represents $odd$ prime numbers $p>2$ with the property that
the sum of three $odd$ primes result in every $odd$ number $n$ greater than 7, $n>7$. ($Odd$ primes $p>2$ may be used more than once in the same sum).


\end{proposition}

 

\begin{proof} 
For {\bf Case I}. $P(x)=N(x)_{n>1}$

\begin{equation}
\label{proof1}
 \left\{{p} + sin \left(\frac{1}{p}\right)\right\} = \left\{{n_{odd_{n>1}}} + sin \left(\frac{1}{n_{n>1}}\right)\right\}
\end{equation}

if the set of all elements are equal, then single elements are equal

\begin{equation}
\label{proof11}
 {p} + sin \left(\frac{1}{p}\right) = {n_{n>1}} + sin \left(\frac{1}{n_{n>1}}\right)
\end{equation}

and the equation can be reduced to 

\begin{equation}
\label{proof111}
 p = n_{n>1}
\end{equation}

by taking out of consideration the periods continuation.

While for {\bf Case II}. $P(x) \iff N(x)_{n>1}$

\begin{equation}
\label{proof2}
 \left\{{p} + sin \left(\frac{1}{p}\right)\right\} \iff \left\{{n_{n>1}} + sin \left(\frac{1}{n_{n>1}}\right)\right\}
\end{equation}

if the set of all elements are equivalent and isomorphic, then single elements are

\begin{equation}
\label{proof22}
 {p} + sin \left(\frac{1}{p}\right) \iff {n_{n>1}} + sin \left(\frac{1}{n_{n>1}}\right)
\end{equation}

and the equation can be reduced to 

\begin{equation}
\label{proof222}
 p \iff n_{n>1}
\end{equation}

 

In the same way, for {\bf Case I}. 
$P_{odd}(x)_{p > 2}=N_{odd}(x)_{n>1}$

\begin{equation}
\label{proof3}
 \left\{ {p_{p>2}} + sin \left(\frac{1}{2p_{p>2}}\right) \right\} = \left\{{n_{odd_{n>1}}} + sin \left(\frac{1}{2n_{odd_{n>1}}}\right)\right\}
\end{equation}

if the set of all elements are equal, then single elements are equal

\begin{equation}
\label{proof33}
  {p_{p>2}} + sin \left(\frac{1}{2p_{p>2}}\right) = {n_{odd_{n>1}}} + sin \left(\frac{1}{2n_{odd_{n>1}}}\right)
\end{equation}

and the equation can be reduced to 

\begin{equation}
\label{proof333}
 p_{p>2} = n_{odd_{n>1}}
\end{equation}

by taking out of consideration the periods continuation.

While for {\bf Case II}. $P_{odd}(x)_{p > 2} \iff N_{odd}(x)_{n>1}$

\begin{equation}
\label{proof4}
 \left\{ {p_{p>2}} + sin \left(\frac{1}{2p_{p>2}}\right) \right\} \iff \left\{{n_{odd_{n>1}}} + sin \left(\frac{1}{2n_{odd_{n>1}}}\right)\right\}
\end{equation}

if the set of all elements are equivalent and isomorphic, then single elements are

\begin{equation}
\label{proof44}
  {p_{p>2}} + sin \left(\frac{1}{2p_{p>2}}\right) \iff {n_{odd_{n>1}}} + sin \left(\frac{1}{2n_{odd_{n>1}}}\right)
\end{equation}

and the equation can be reduced to 

\begin{equation}
\label{proof444}
 p_{p>2} \iff n_{odd_{n>1}}
\end{equation}

Likewise, the weak Goldbach's conjecture is proven, 
by been able to zero-cross all $odd$ natural numbers greater than seven, $n_{odd}> 7$, starting from $odd$ primes $p>2$, as with $odd$ numbers greater than 1, $n_{odd}> 1$.

\end{proof}



\begin{thebibliography}{9}
\bibitem{hardy_cite}
Hardy G. H. and Littlewood J. E., Some problems of ÒPartitio NumerorumÓ III: On the expression of a number as a sum of primes, Acta. Math., 44 (1923). pp.1-70.

\bibitem{gove_cite}
Effinger G., Some numerical implication of the Hardy and Littlewood analysis of the 3- primes problem, submitted for publication.

\bibitem{bruno_cite}
Lucke B., Zur Hardy-Littlewoodschen Behandlung des Goldbachschen Problems,Ó Doctoral Dissertation, Gottingen, 1926.

\bibitem{dimi_cite}
Zinoviev D., On Vinogradov's constant in Goldbach's ternary problem, Journal of Number Theory, 65 (1997), 334Ð358.

\bibitem{boro_cite}
Borodzkin K. G., On I. M. VinogradoiÞs constant, Proc. 3rd All-Union Math. Conf., vol. 1,. Izdat. Akad. Nauk SSSR, Moscow, 1956. (Russian) MR 20:6973a.

\bibitem{chen_cite}
Chen J. R.  and Wang T. Z., On the odd Goldbach problem, Acta Math. Sinica 32 (1989), 702-718.

\bibitem{chen2_cite}
Chen J. R.  and Wang T. Z., The Goldbach problem for odd numbers, Acta Math. Sinica 39 (1996), 169-174.

\bibitem{x1_cite}
Deshouillers J. M., te Riele H. J. and Saouter Y., New Experimental Results Concerning the Goldbach Conjecture, In ÒProc. 3rd Int. Symp. on Algorithm Number TheoryÓ, LNCS, 1423 (1998), 204-215.

\bibitem{x2_cite}
Richstein J., Verifying GoldbachÕs conjecture up to $4 \times 10^{14}$ , Mathematics and Computation 70 (2001), 1745-1749.

\bibitem{tao_cite}
Tao T., Every Odd number Greater Than 1 is the Sum of at most five Primes, 
\url{http://arXiv:1201.6656v3}

\bibitem{ramare_cite}
Ramar\'{e} O. and Saouter Y., Short effective intervals containing primes, J. Number Theory 98 (2003), no. 1, 10-33.

\bibitem{eulerproduct1_cite}
Edwards H. M., 
The Euler Product Formula,  
New York: Dover 2001, pp. 6-7.

\bibitem{eulerproduct2_cite}
Shimura G., 
Euler Products and Eisenstein Series, 
Providence, RI: Amer. Math. Soc., 1997.

\bibitem{sieves} 
O'Neill M. E., 
The Genuine Sieve of Eratosthenes, 
Journal of Functional Programming, 
Cambridge University Press 2008,
pp. 10-11, 
DOI:10.1017/S0956796808007004.

\bibitem{dynamicaleras}
Mateos L.A., Dynamical Sieve of Eratosthenes, 
\url{http://arxiv.org/pdf/1206.2791}

\bibitem{mateos1_cite} 
Mateos L.A., Chaotic Nonlinear Prime Number Function, 
AIP Conf. Proc. CMLS 2011, 1371, pp. 161-170, 
DOI:10.1063/1.3596639.

\bibitem{eras_cite} 
Sieve Of Eratosthenes. \url{http://c2.com/cgi/wiki?SieveOfEratosthenes}. \\
Retrieved on 2012-02-15.


\end{thebibliography}
\end{document}